\documentclass[12pt,reqno,a4wide]{amsart}

\usepackage{tikz} 
\usepackage{calrsfs}
\usepackage{mathrsfs}

\usepackage{array}
\usepackage{hyperref}

\usetikzlibrary{decorations.pathreplacing}
\usetikzlibrary{fit}

\usepackage{pgfplots}

\allowdisplaybreaks

\catcode`,\active

\catcode`\,12

\begin{document}
 \bibliographystyle{plain}
 
 %
 %
 
 \title
  {The height of Dyck paths and checkerboard labellings}

 \author[H. Prodinger ]{Helmut Prodinger }
 \address{Department of Mathematics, University of Stellenbosch 7602, Stellenbosch, South Africa
  and
  NITheCS (National Institute for
  Theoretical and Computational Sciences), South Africa.}
 \email{warrenham33@gmail.com}


 \begin{abstract}
Dyck paths and certain black/white labelling of nodes leads to the \emph{white-height}. Using generating functions and tricks of the trade, 
we establish that the average white-height among Dyck paths of half length $n$ is asymptotic to $\frac12\sqrt{\pi n}$ for two different models. 
These are appealing results that could be presented to students to learn the trade. A last section links walks with double steps and 2-Motzkin paths. For them, the white-height
is the natural concept.

Folks who might be offended by the notion of \emph{white-height} might choose their own colours that they like.
 \end{abstract}.

 \maketitle

\section{Introduction}

Dyck paths and Dyck numbers are among the most popular combinatorial objects, see \cite{OEIS} and \cite{stanley-book}.
Recently I have seen some papers where a black/white pattern is combined with Dyck paths; we give just two citations
\cite{Baril} and \cite{china}. I will consider the \emph{white}-height of such Dyck paths, which is the maximum of white nodes under each node of the Dyck path.
Ignoring the black/white discussion, this will be the \emph{height} ($+1$ to be exact), a subject that was discussed in the pioneering paper \cite{deBrKnRi}.
We will show how to do exact enumeration of Dyck paths with white-height $\ge h$ (or $= h$, $\le h$). 

We like to use the term \emph{checkerboard} and consider two models. In the first one, the parity of the levels decides the colour. To get the checkerboard feel, we must rotate the Dyck path by 45 degrees. In the second model, the colour is decided by the parity of abscissa $+$ ordinate, leading immediately to a checkerboard scheme. 

The enumeration (appropriate generating functions) follows the classic paper \cite{deBrKnRi} but with some interesting twists that the reader will like. 
To compute the average white-height, we will use a combined scheme of Mellin transform and singularity analysis, that I used already about 45 years ago \cite{prodinger-ars}; however, we like to cite a more recent publication \cite{HPW} as it discussed some technical  details (of a related but different problem) at some length.

We are optimistic that this presentation serves a pedagogic value as well. 

\section{Labels based on parity of ordinates}
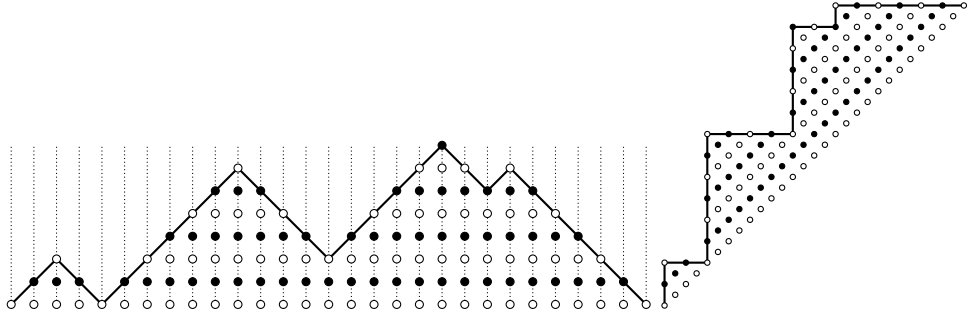
\begin{figure}[h]
 \begin{tikzpicture}[scale=0.3,rotate=0]
  \foreach\x in {0,...,28}
  {\draw[fill=white][ thin, densely dotted](\x,0)--(\x,7);}

  \draw[fill=white] [thick] (0,0) -- (1,1)--(2,2);
  
  \draw[fill=white] [thick] (2,2)--(3,1);
  \draw[fill=white] [thick] (3,1)--(4,0);
  \draw[fill=white] [thick] (4,0)--(10,6);
  \draw[fill=white] [thick] (10,6)--(11,5);
  \draw[fill=white] [thick] (11,5)--(13,3);
  \draw[fill=white] [thick ] (13,3)--(14,2);
  \draw[fill=white] [thick] (14,2)--(19,7);
  \draw[fill=white] [thick ] (19,7)--(21,5);
  \draw[fill=white] [thick] (21,5)--(22,6);
  \draw[fill=white] [thick ] (22,6)--(24,4);
  \draw[fill=white] [thick ] (24,4)--(27,1);
  \draw[fill=white] [thick ] (27,1)--(28,0);

  \foreach \x in {0}
  \foreach \y in {0}
  {
   \draw[fill=white] (\x,\y) circle (5pt);
  }
  
  \foreach \x in {1}   \foreach \y in {1}{    \filldraw (\x,\y) circle (5pt);   }
  \foreach \x in {1}   \foreach \y in {0}{ \draw[fill=white]    (\x,\y) circle (5pt);   }
  
  \foreach \x in {2}   \foreach \y in {0,2}{    \draw[fill=white] (\x,\y) circle (5pt);   }
  \foreach \x in {2}   \foreach \y in {1}{   \filldraw  (\x,\y) circle (5pt);   }
  
  \foreach \x in {3}   \foreach \y in {1}{    \filldraw (\x,\y) circle (5pt);   }
  \foreach \x in {3}   \foreach \y in {0}{   \draw[fill=white]  (\x,\y) circle (5pt);   }

  \foreach \x in {4}   \foreach \y in {0}{    \draw[fill=white] (\x,\y) circle (5pt);   }    
  \foreach \x in {5}   \foreach \y in {1}{    \filldraw (\x,\y) circle (5pt);   } 
  \foreach \x in {5}   \foreach \y in {0}{   \draw[fill=white]  (\x,\y) circle (5pt);   }    
  \foreach \x in {6}   \foreach \y in {0,2}{    \draw[fill=white] (\x,\y) circle (5pt);   } 
  \foreach \x in {6}   \foreach \y in {1}{   \filldraw  (\x,\y) circle (5pt);   }    
  \foreach \x in {7}   \foreach \y in {1,3}{    \filldraw (\x,\y) circle (5pt);   }   
  \foreach \x in {7}   \foreach \y in {0,2}{    \draw[fill=white] (\x,\y) circle (5pt);   }   
  \foreach \x in {8}   \foreach \y in {0,2,4}{    \draw[fill=white] (\x,\y) circle (5pt);   }     
  \foreach \x in {8}   \foreach \y in {1,3}{    \filldraw (\x,\y) circle (5pt);   }     
  \foreach \x in {9}   \foreach \y in {1,3,5}{    \filldraw (\x,\y) circle (5pt);   }     
  \foreach \x in {9}   \foreach \y in {0,2,4}{    \draw[fill=white] (\x,\y) circle (5pt);   }     
  \foreach \x in {10}   \foreach \y in {1,3,5}{    \filldraw (\x,\y) circle (5pt);   }     
  \foreach \x in {10}   \foreach \y in {0,2,4,6}{    \draw[fill=white] (\x,\y) circle (5pt);   }     
  \foreach \x in {11}   \foreach \y in {1,3,5}{    \filldraw (\x,\y) circle (5pt);   }     
  \foreach \x in {11}   \foreach \y in {0,2,4}{    \draw[fill=white] (\x,\y) circle (5pt);   } 
  \foreach \x in {12}   \foreach \y in {1,3}{    \filldraw (\x,\y) circle (5pt);   }     
  \foreach \x in {12}   \foreach \y in {0,2,4}{    \draw[fill=white] (\x,\y) circle (5pt);   }
  \foreach \x in {13}   \foreach \y in {1,3}{    \filldraw (\x,\y) circle (5pt);   }     
  \foreach \x in {13}   \foreach \y in {0,2}{    \draw[fill=white] (\x,\y) circle (5pt);   } 
  \foreach \x in {14}   \foreach \y in {1}{    \filldraw (\x,\y) circle (5pt);   }     
  \foreach \x in {14}   \foreach \y in {0,2}{    \draw[fill=white] (\x,\y) circle (5pt);   }     
  \foreach \x in {15}   \foreach \y in {1,3}{    \filldraw (\x,\y) circle (5pt);   }     
  \foreach \x in {15}   \foreach \y in {0,2}{    \draw[fill=white] (\x,\y) circle (5pt);   }
  \foreach \x in {16}   \foreach \y in {1,3}{    \filldraw (\x,\y) circle (5pt);   }     
  \foreach \x in {16}   \foreach \y in {0,2,4}{    \draw[fill=white] (\x,\y) circle (5pt);   } 
  \foreach \x in {17}   \foreach \y in {1,3,5}{    \filldraw (\x,\y) circle (5pt);   }     
  \foreach \x in {17}   \foreach \y in {0,2,4}{    \draw[fill=white] (\x,\y) circle (5pt);   } 
  \foreach \x in {18}   \foreach \y in {1,3,5}{    \filldraw (\x,\y) circle (5pt);   }     
  \foreach \x in {18}   \foreach \y in {0,2,4,6}{    \draw[fill=white] (\x,\y) circle (5pt);   } 
  \foreach \x in {19}   \foreach \y in {1,3,5,7}{    \filldraw (\x,\y) circle (5pt);   }     
  \foreach \x in {19}   \foreach \y in {0,2,4,6}{    \draw[fill=white] (\x,\y) circle (5pt);   } 
  \foreach \x in {20}   \foreach \y in {1,3,5}{    \filldraw (\x,\y) circle (5pt);   }     
  \foreach \x in {20}   \foreach \y in {0,2,4,6}{    \draw[fill=white] (\x,\y) circle (5pt);   }   
  \foreach \x in {21}   \foreach \y in {1,3,5}{    \filldraw (\x,\y) circle (5pt);   }     
  \foreach \x in {21}   \foreach \y in {0,2,4}{    \draw[fill=white] (\x,\y) circle (5pt);   } 
  \foreach \x in {22}   \foreach \y in {1,3,5}{    \filldraw (\x,\y) circle (5pt);   }     
  \foreach \x in {22}   \foreach \y in {0,2,4,6}{    \draw[fill=white] (\x,\y) circle (5pt);   } 
  \foreach \x in {23}   \foreach \y in {1,3,5}{    \filldraw (\x,\y) circle (5pt);   }     
  \foreach \x in {23}   \foreach \y in {0,2,4}{    \draw[fill=white] (\x,\y) circle (5pt);   } 
  
  \foreach \x in {24}   \foreach \y in {1,3}{    \filldraw (\x,\y) circle (5pt);   }     
  \foreach \x in {24}   \foreach \y in {0,2,4}{    \draw[fill=white] (\x,\y) circle (5pt);   } 
  \foreach \x in {25}   \foreach \y in {1,3}{    \filldraw (\x,\y) circle (5pt);   }     
  \foreach \x in {25}   \foreach \y in {0,2}{    \draw[fill=white] (\x,\y) circle (5pt);   } 
  
  \foreach \x in {26}   \foreach \y in {1}{    \filldraw (\x,\y) circle (5pt);   }     
  \foreach \x in {26}   \foreach \y in {0,2}{    \draw[fill=white] (\x,\y) circle (5pt);   } 
  \foreach \x in {27}   \foreach \y in {1}{    \filldraw (\x,\y) circle (5pt);   }     
  \foreach \x in {27}   \foreach \y in {0}{    \draw[fill=white] (\x,\y) circle (5pt);   } 
  
  \foreach \x in {28}   \foreach \y in {}{    \filldraw (\x,\y) circle (5pt);   }     
  \foreach \x in {28}   \foreach \y in {0}{    \draw[fill=white] (\x,\y) circle (5pt);   } 
 \end{tikzpicture}
 \begin{tikzpicture}[scale=0.2,rotate=45]

 \draw[fill=white] [thick] (0,0) -- (1,1)--(2,2);
 
 \draw[fill=white] [thick] (2,2)--(3,1);
 \draw[fill=white] [thick] (3,1)--(4,0);
 \draw[fill=white] [thick] (4,0)--(10,6);
 \draw[fill=white] [thick] (10,6)--(11,5);
 \draw[fill=white] [thick] (11,5)--(13,3);
 \draw[fill=white] [thick ] (13,3)--(14,2);
 \draw[fill=white] [thick] (14,2)--(19,7);
 \draw[fill=white] [thick ] (19,7)--(21,5);
 \draw[fill=white] [thick] (21,5)--(22,6);
 \draw[fill=white] [thick ] (22,6)--(24,4);
 \draw[fill=white] [thick ] (24,4)--(27,1);
 \draw[fill=white] [thick ] (27,1)--(28,0);

 \foreach \x in {0}
 \foreach \y in {0}
 {
  \draw[fill=white] (\x,\y) circle (5pt);
 }
 
 \foreach \x in {1}   \foreach \y in {1}{    \filldraw (\x,\y) circle (5pt);   }
 \foreach \x in {1}   \foreach \y in {0}{ \draw[fill=white]    (\x,\y) circle (5pt);   }
 
 \foreach \x in {2}   \foreach \y in {0,2}{    \draw[fill=white] (\x,\y) circle (5pt);   }
 \foreach \x in {2}   \foreach \y in {1}{   \filldraw  (\x,\y) circle (5pt);   }
 
 \foreach \x in {3}   \foreach \y in {1}{    \filldraw (\x,\y) circle (5pt);   }
 \foreach \x in {3}   \foreach \y in {0}{   \draw[fill=white]  (\x,\y) circle (5pt);   }

 \foreach \x in {4}   \foreach \y in {0}{    \draw[fill=white] (\x,\y) circle (5pt);   }    
 \foreach \x in {5}   \foreach \y in {1}{    \filldraw (\x,\y) circle (5pt);   } 
 \foreach \x in {5}   \foreach \y in {0}{   \draw[fill=white]  (\x,\y) circle (5pt);   }    
 \foreach \x in {6}   \foreach \y in {0,2}{    \draw[fill=white] (\x,\y) circle (5pt);   } 
 \foreach \x in {6}   \foreach \y in {1}{   \filldraw  (\x,\y) circle (5pt);   }    
 \foreach \x in {7}   \foreach \y in {1,3}{    \filldraw (\x,\y) circle (5pt);   }   
 \foreach \x in {7}   \foreach \y in {0,2}{    \draw[fill=white] (\x,\y) circle (5pt);   }   
 \foreach \x in {8}   \foreach \y in {0,2,4}{    \draw[fill=white] (\x,\y) circle (5pt);   }     
 \foreach \x in {8}   \foreach \y in {1,3}{    \filldraw (\x,\y) circle (5pt);   }     
 \foreach \x in {9}   \foreach \y in {1,3,5}{    \filldraw (\x,\y) circle (5pt);   }     
 \foreach \x in {9}   \foreach \y in {0,2,4}{    \draw[fill=white] (\x,\y) circle (5pt);   }     
 \foreach \x in {10}   \foreach \y in {1,3,5}{    \filldraw (\x,\y) circle (5pt);   }     
 \foreach \x in {10}   \foreach \y in {0,2,4,6}{    \draw[fill=white] (\x,\y) circle (5pt);   }     
 \foreach \x in {11}   \foreach \y in {1,3,5}{    \filldraw (\x,\y) circle (5pt);   }     
 \foreach \x in {11}   \foreach \y in {0,2,4}{    \draw[fill=white] (\x,\y) circle (5pt);   } 
 \foreach \x in {12}   \foreach \y in {1,3}{    \filldraw (\x,\y) circle (5pt);   }     
 \foreach \x in {12}   \foreach \y in {0,2,4}{    \draw[fill=white] (\x,\y) circle (5pt);   }
 \foreach \x in {13}   \foreach \y in {1,3}{    \filldraw (\x,\y) circle (5pt);   }     
 \foreach \x in {13}   \foreach \y in {0,2}{    \draw[fill=white] (\x,\y) circle (5pt);   } 
 \foreach \x in {14}   \foreach \y in {1}{    \filldraw (\x,\y) circle (5pt);   }     
 \foreach \x in {14}   \foreach \y in {0,2}{    \draw[fill=white] (\x,\y) circle (5pt);   }     
 \foreach \x in {15}   \foreach \y in {1,3}{    \filldraw (\x,\y) circle (5pt);   }     
 \foreach \x in {15}   \foreach \y in {0,2}{    \draw[fill=white] (\x,\y) circle (5pt);   }
 \foreach \x in {16}   \foreach \y in {1,3}{    \filldraw (\x,\y) circle (5pt);   }     
 \foreach \x in {16}   \foreach \y in {0,2,4}{    \draw[fill=white] (\x,\y) circle (5pt);   } 
 \foreach \x in {17}   \foreach \y in {1,3,5}{    \filldraw (\x,\y) circle (5pt);   }     
 \foreach \x in {17}   \foreach \y in {0,2,4}{    \draw[fill=white] (\x,\y) circle (5pt);   } 
 \foreach \x in {18}   \foreach \y in {1,3,5}{    \filldraw (\x,\y) circle (5pt);   }     
 \foreach \x in {18}   \foreach \y in {0,2,4,6}{    \draw[fill=white] (\x,\y) circle (5pt);   } 
 \foreach \x in {19}   \foreach \y in {1,3,5,7}{    \filldraw (\x,\y) circle (5pt);   }     
 \foreach \x in {19}   \foreach \y in {0,2,4,6}{    \draw[fill=white] (\x,\y) circle (5pt);   } 
 \foreach \x in {20}   \foreach \y in {1,3,5}{    \filldraw (\x,\y) circle (5pt);   }     
 \foreach \x in {20}   \foreach \y in {0,2,4,6}{    \draw[fill=white] (\x,\y) circle (5pt);   }   
 \foreach \x in {21}   \foreach \y in {1,3,5}{    \filldraw (\x,\y) circle (5pt);   }     
 \foreach \x in {21}   \foreach \y in {0,2,4}{    \draw[fill=white] (\x,\y) circle (5pt);   } 
 \foreach \x in {22}   \foreach \y in {1,3,5}{    \filldraw (\x,\y) circle (5pt);   }     
 \foreach \x in {22}   \foreach \y in {0,2,4,6}{    \draw[fill=white] (\x,\y) circle (5pt);   } 
 \foreach \x in {23}   \foreach \y in {1,3,5}{    \filldraw (\x,\y) circle (5pt);   }     
 \foreach \x in {23}   \foreach \y in {0,2,4}{    \draw[fill=white] (\x,\y) circle (5pt);   } 
 
 \foreach \x in {24}   \foreach \y in {1,3}{    \filldraw (\x,\y) circle (5pt);   }     
 \foreach \x in {24}   \foreach \y in {0,2,4}{    \draw[fill=white] (\x,\y) circle (5pt);   } 
 \foreach \x in {25}   \foreach \y in {1,3}{    \filldraw (\x,\y) circle (5pt);   }     
 \foreach \x in {25}   \foreach \y in {0,2}{    \draw[fill=white] (\x,\y) circle (5pt);   } 
 
 \foreach \x in {26}   \foreach \y in {1}{    \filldraw (\x,\y) circle (5pt);   }     
 \foreach \x in {26}   \foreach \y in {0,2}{    \draw[fill=white] (\x,\y) circle (5pt);   } 
 \foreach \x in {27}   \foreach \y in {1}{    \filldraw (\x,\y) circle (5pt);   }     
 \foreach \x in {27}   \foreach \y in {0}{    \draw[fill=white] (\x,\y) circle (5pt);   } 
 
 \foreach \x in {28}   \foreach \y in {}{    \filldraw (\x,\y) circle (5pt);   }     
 \foreach \x in {28}   \foreach \y in {0}{    \draw[fill=white] (\x,\y) circle (5pt);   } 
\end{tikzpicture}
 
 \caption{A Dyck path with black and white nodes. The same, but rotated by 45 degrees, so that it looks more `checkerboard.'}\label{susi1}
\end{figure}

In Figure~\ref{susi1} we see a Dyck path where the nodes between the path and the $x$-axis are label black/white as indicated. To get more
of a `checkerboard' feeling is was redrawn after a rotation by 45 degrees.

The \emph{white-height} of an $x$-coordinate is the number of white nodes on that coordinate. The \emph{white-height} of the whole Dyck path is the maximum of all these values. Since the white-height is more or less half of the ordinary height, we expect an asymptotic value $\sim\frac12\sqrt{\pi n}$ when all Dyck paths of half-length $n$ are equally likely,  but we will provide more precise calculations.

A typical decomposition of a Dyck path is to consider the \emph{first return} to the $x$-axis. Then the first and the last step of the first sojurn can be chopped off, leading to a recursion. But careful!
This leads then to a Dyck path where the labels are black-white--\dots instead of the original white-black--\dots\;. So we introduce two sequences of generating functions:
$A_h(z)$ counts paths white-black--\dots paths with \emph{white-height} $\le h$ and
$B_h(z)$ counts paths black-white--\dots paths with \emph{white-height} $\le h$. These functions call each other recursively:
\begin{align*}
 A_h&=1+zB_{h-1}A_h,\quad A_0=0,\\
 B_h&=1+zA_hB_h,\quad B_0=1
\end{align*}
or
\begin{align*}
 A_h&=\frac{1}{1-zB_{h-1}},\quad A_0=0,\\
 B_h&=\frac{1}{1-zA_{h}},\quad B_0=1.
\end{align*}
One could now eliminate the $B_h$ sequence to get
\begin{equation*}
A_h=\cfrac{1}{1-\cfrac{z}{1-zA_{h-1}}},
\end{equation*}
which has the pleasant continued fraction appeal.

After doing a few calculations, one sees the quantity $\sqrt{1-4z}$ popping up, and so it is time to introduce the substitution
\begin{equation*}
 z=\frac{u}{(1+u)^2}
\end{equation*}
that is successfully known since 1972 \cite{deBrKnRi} (and sometimes rediscovered). Then we find
\begin{align*}
 A_h&=(1+u)\frac{1-u^{2h}}{1-u^{2h+1}},\\
 B_h&=(1+u)\frac{1-u^{2h+1}}{1-u^{2h+2}}.
\end{align*}
People might be surprised how these were found. Numerators and denominators satisfy second order recursions in all instances and can be guessed/solved. For such guessing approaches we use the program \textsf{gfun}.
Having the final forms, they may be proved by induction in a routinely fashion.

We are only interested in $A_h(z)$ and consider $B_h(z)$ just as an auxiliary variable.
Pushing $h$ to infinity means that there is no restriction w.r.t. the white-height, and so all Dyck paths should be counted. Indeed,
$A_\infty=1+u=\frac{1-\sqrt{1-4z}}{2z}$. The difference $A_\infty-A_h$ refers to all Dyck paths with white-height $>h$, a quantity that comes handy
when one wants to compute the average white-height.
\begin{equation*}
 A_\infty-A_{h}=(1+u)\frac{1-u^{2h+1}-1+u^{2h}}{1-u^{2h+1}}=(1+u)\frac{u^{2h}(1-u)}{1-u^{2h+1}}=\frac{1-u^2}{u}\frac{u^{2h+1}}{1-u^{2h+1}}.
\end{equation*}
 
\begin{figure}[h]
 \begin{tikzpicture}[scale=0.3]
  \foreach\x in {0,...,28}
  {\draw[fill=white][ thin, densely dotted](\x,0)--(\x,7);}

  \draw[fill=white] [thick] (0,0) -- (1,1)--(2,2);
  
  \draw[fill=white] [thick] (2,2)--(3,1);
  \draw[fill=white] [thick] (3,1)--(4,0);
  \draw[fill=white] [thick] (4,0)--(10,6);
  \draw[fill=white] [thick] (10,6)--(11,5);
  \draw[fill=white] [thick] (11,5)--(13,3);
  \draw[fill=white] [thick ] (13,3)--(14,2);
  \draw[fill=white] [thick] (14,2)--(19,7);
  \draw[fill=white] [thick ] (19,7)--(21,5);
  \draw[fill=white] [thick] (21,5)--(22,6);
  \draw[fill=white] [thick ] (22,6)--(24,4);
  \draw[fill=white] [thick ] (24,4)--(27,1);
  \draw[fill=white] [thick ] (27,1)--(28,0);

  \foreach \x in {0}
  \foreach \y in {0}
  {
   \filldraw (\x,\y) circle (5pt);
  }
  
  \foreach \x in {1}   \foreach \y in {1}{    \draw[fill=white] (\x,\y) circle (5pt);   }
  \foreach \x in {1}   \foreach \y in {0}{ \filldraw    (\x,\y) circle (5pt);   }
  
  \foreach \x in {2}   \foreach \y in {0,2}{    \filldraw (\x,\y) circle (5pt);   }
  \foreach \x in {2}   \foreach \y in {1}{   \draw[fill=white]  (\x,\y) circle (5pt);   }
  
  \foreach \x in {3}   \foreach \y in {1}{    \draw[fill=white] (\x,\y) circle (5pt);   }
  \foreach \x in {3}   \foreach \y in {0}{   \filldraw  (\x,\y) circle (5pt);   }

  \foreach \x in {4}   \foreach \y in {0}{    \filldraw (\x,\y) circle (5pt);   }    
  \foreach \x in {5}   \foreach \y in {1}{    \draw[fill=white] (\x,\y) circle (5pt);   } 
  \foreach \x in {5}   \foreach \y in {0}{   \filldraw  (\x,\y) circle (5pt);   }    
  \foreach \x in {6}   \foreach \y in {0,2}{    \filldraw (\x,\y) circle (5pt);   } 
  \foreach \x in {6}   \foreach \y in {1}{   \draw[fill=white]  (\x,\y) circle (5pt);   }    
  \foreach \x in {7}   \foreach \y in {1,3}{    \draw[fill=white] (\x,\y) circle (5pt);   }   
  \foreach \x in {7}   \foreach \y in {0,2}{    \filldraw (\x,\y) circle (5pt);   }   
  \foreach \x in {8}   \foreach \y in {0,2,4}{    \filldraw (\x,\y) circle (5pt);   }     
  \foreach \x in {8}   \foreach \y in {1,3}{    \draw[fill=white] (\x,\y) circle (5pt);   }     
  \foreach \x in {9}   \foreach \y in {1,3,5}{    \draw[fill=white] (\x,\y) circle (5pt);   }     
  \foreach \x in {9}   \foreach \y in {0,2,4}{    \filldraw (\x,\y) circle (5pt);   }     
  \foreach \x in {10}   \foreach \y in {1,3,5}{    \draw[fill=white] (\x,\y) circle (5pt);   }     
  \foreach \x in {10}   \foreach \y in {0,2,4,6}{    \filldraw (\x,\y) circle (5pt);   }     
  \foreach \x in {11}   \foreach \y in {1,3,5}{    \draw[fill=white] (\x,\y) circle (5pt);   }     
  \foreach \x in {11}   \foreach \y in {0,2,4}{    \filldraw (\x,\y) circle (5pt);   } 
  \foreach \x in {12}   \foreach \y in {1,3}{    \draw[fill=white] (\x,\y) circle (5pt);   }     
  \foreach \x in {12}   \foreach \y in {0,2,4}{    \filldraw (\x,\y) circle (5pt);   }
  \foreach \x in {13}   \foreach \y in {1,3}{    \draw[fill=white] (\x,\y) circle (5pt);   }     
  \foreach \x in {13}   \foreach \y in {0,2}{    \filldraw (\x,\y) circle (5pt);   } 
  \foreach \x in {14}   \foreach \y in {1}{    \draw[fill=white] (\x,\y) circle (5pt);   }     
  \foreach \x in {14}   \foreach \y in {0,2}{    \filldraw (\x,\y) circle (5pt);   }     
  \foreach \x in {15}   \foreach \y in {1,3}{    \draw[fill=white] (\x,\y) circle (5pt);   }     
  \foreach \x in {15}   \foreach \y in {0,2}{    \filldraw (\x,\y) circle (5pt);   }
  \foreach \x in {16}   \foreach \y in {1,3}{    \draw[fill=white] (\x,\y) circle (5pt);   }     
  \foreach \x in {16}   \foreach \y in {0,2,4}{    \filldraw (\x,\y) circle (5pt);   } 
  \foreach \x in {17}   \foreach \y in {1,3,5}{    \draw[fill=white] (\x,\y) circle (5pt);   }     
  \foreach \x in {17}   \foreach \y in {0,2,4}{    \filldraw (\x,\y) circle (5pt);   } 
  \foreach \x in {18}   \foreach \y in {1,3,5}{    \draw[fill=white] (\x,\y) circle (5pt);   }     
  \foreach \x in {18}   \foreach \y in {0,2,4,6}{    \filldraw (\x,\y) circle (5pt);   } 
  \foreach \x in {19}   \foreach \y in {1,3,5,7}{    \draw[fill=white] (\x,\y) circle (5pt);   }     
  \foreach \x in {19}   \foreach \y in {0,2,4,6}{    \filldraw (\x,\y) circle (5pt);   } 
  \foreach \x in {20}   \foreach \y in {1,3,5}{    \draw[fill=white] (\x,\y) circle (5pt);   }     
  \foreach \x in {20}   \foreach \y in {0,2,4,6}{    \filldraw (\x,\y) circle (5pt);   }   
  \foreach \x in {21}   \foreach \y in {1,3,5}{    \draw[fill=white] (\x,\y) circle (5pt);   }     
  \foreach \x in {21}   \foreach \y in {0,2,4}{    \filldraw (\x,\y) circle (5pt);   } 
  \foreach \x in {22}   \foreach \y in {1,3,5}{    \draw[fill=white] (\x,\y) circle (5pt);   }     
  \foreach \x in {22}   \foreach \y in {0,2,4,6}{    \filldraw (\x,\y) circle (5pt);   } 
  \foreach \x in {23}   \foreach \y in {1,3,5}{    \filldraw (\x,\y) circle (5pt);   }     
  \foreach \x in {23}   \foreach \y in {0,2,4}{    \filldraw (\x,\y) circle (5pt);   } 
  
  \foreach \x in {24}   \foreach \y in {1,3}{    \draw[fill=white] (\x,\y) circle (5pt);   }     
  \foreach \x in {24}   \foreach \y in {0,2,4}{    \filldraw (\x,\y) circle (5pt);   } 
  \foreach \x in {25}   \foreach \y in {1,3}{    \draw[fill=white] (\x,\y) circle (5pt);   }     
  \foreach \x in {25}   \foreach \y in {0,2}{    \filldraw (\x,\y) circle (5pt);   } 
  
  \foreach \x in {26}   \foreach \y in {1}{    \draw[fill=white] (\x,\y) circle (5pt);   }     
  \foreach \x in {26}   \foreach \y in {0,2}{    \filldraw (\x,\y) circle (5pt);   } 
  \foreach \x in {27}   \foreach \y in {1}{    \draw[fill=white] (\x,\y) circle (5pt);   }     
  \foreach \x in {27}   \foreach \y in {0}{    \filldraw (\x,\y) circle (5pt);   } 
  
  \foreach \x in {28}   \foreach \y in {}{    \draw[fill=white] (\x,\y) circle (5pt);   }     
  \foreach \x in {28}   \foreach \y in {0}{    \filldraw (\x,\y) circle (5pt);   } 
 \end{tikzpicture}
 
 \caption{A path related to the functions $B_h(z)$.}
\end{figure}
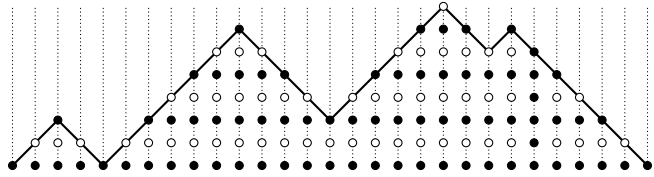

We will show how to extract the coefficients of $A_\infty-A_h$. This follows the method that was used in \cite{deBrKnRi}.
\begin{align*}
[z^n](A_\infty-A_h)&=[z^n](A_\infty-A_h)\\
&=\frac1{2\pi i}\oint\frac{dz}{z^{n+1}}\frac{1-u^2}{u}\frac{u^{2h+1}}{1-u^{2h+1}}\\
&=\frac1{2\pi i}\oint\frac{du(1-u)(1+u)^{2n+2}}{(1+u)^3u^{n+1}}\frac{1-u^2}{u}\sum_{k\ge1}u^{k(2h+1)}\\
&=\sum_{k\ge1}\frac1{2\pi i}\oint\frac{du(1-u)^2(1+u)^{2n}}{u^{n+2-k(2h+1)}}\\
&=\sum_{k\ge1}[u^{n+1-k(2h+1)}](1-2u+u^2)(1+u)^{2n}\\
&=\sum_{k\ge1}\biggl[\binom{2n}{n+1-k(2h+1)}-2\binom{2n}{n-k(2h+1)}+\binom{2n}{n-1-k(2h+1)}\biggr].
\end{align*}
For the average, we must sum:
\begin{equation*}
\sum_{h\ge1}\sum_{k\ge1}\biggl[\binom{2n}{n+1-k(2h+1)}-2\binom{2n}{n-k(2h+1)}+\binom{2n}{n-1-k(2h+1)}\biggr].
\end{equation*}
To do asymptotics, it is more convenient to start from
\begin{equation*}
\sum_{h\ge1}\frac{1-u^2}{u}\frac{u^{2h+1}}{1-u^{2h+1}}=
\frac{1-u^2}{u}\sum_{h\ge0}\frac{u^{2h+1}}{1-u^{2h+1}}-(1+u).
\end{equation*}
One sets $u=e^{-t}$ and studies
\begin{equation*}
\sum_{h\ge0}\frac{u^{2h+1}}{1-u^{2h+1}}=\sum_{h\ge0}\sum_{k\ge1}e^{-k(2h+1)t}.
\end{equation*}
The next step is to consider/compute the celebrated Mellin transform. A general reference for Combinatorics with an analytic flavour is
\cite{FS}. The reader might again follow \cite{HPW};
\begin{equation*}
\mathcal{M}\sum_{h\ge0}\sum_{k\ge1}e^{-k(2h+1)t}=\Gamma(s)\sum_{h\ge0}\sum_{k\ge1}k^{-s}(2h+1)^{-s}
=\Gamma(s)\zeta^2(s)(1-2^{-s}).
\end{equation*}
By the inverse Mellin transform we get the function back:
\begin{equation*}
\mathcal {S}:=\sum_{h\ge0}\sum_{k\ge1}e^{-k(2h+1)t}=\frac1{2\pi i}\int_{2-i\infty}^{2+i\infty}\Gamma(s)\zeta^2(s)(1-2^{-s})t^{-s}ds
\end{equation*}
which is approximated by shifting the line of integration to the left and collecting residues. We consider $s=1,0,-1$ and use a computer;
\begin{equation*}
\mathcal {S}\sim\frac{-\log(t)+\log(2)+\gamma}{2t}+\frac{t}{144}.
\end{equation*}
Since 
\begin{equation*}
u=-\log(t)=\frac{1-\sqrt{1-4z}}{2z}-1,
\end{equation*}
we may at the end of the day express our expansions in terms of $\sqrt{1-4z}\,$;
\begin{equation*}
-\log  ( t )+\log  ( 2 ) +\gamma -2+t+ \Bigl(\frac{-\log ( t)}{6} -	\frac {35}{72}+\frac{\log ( 2)+\gamma}{6}   \Bigr) {t}^{2}
\end{equation*}
or, in the $z$-world,
\begin{equation*}
-\frac12\log(1-4z)-2\sqrt{1-4z}-2+\textsf{further terms}.
\end{equation*}
Singularity analysis of generating functions \cite{FO} allows to translate the expansion into an asymptotic expansion of the coefficients. This leads to 
\begin{equation*}
\frac{4^n}{n}+\frac{4^n}{\sqrt{\pi}n^{3/2}}+\textsf{smaller order terms}.
\end{equation*}
However, this must be divided by Catalan numbers
\begin{equation*}
\frac{1}{n+1}\binom{2n}{n}=\frac{4^n}{\sqrt{\pi}n^{3/2}}+\textsf{smaller order terms}
\end{equation*}
After this division one finds the average white-height among all Dyck-paths of semilength $n$  to be asymptotic to
\begin{equation*}
\textsf{Average white-height}\sim	\frac12\sqrt{\pi n}+1+\textsf{smaller order terms}.
\end{equation*}

\section{The checkerboard model}

Now we move to the model where no rotation is necessary to see the checkerboard pattern.

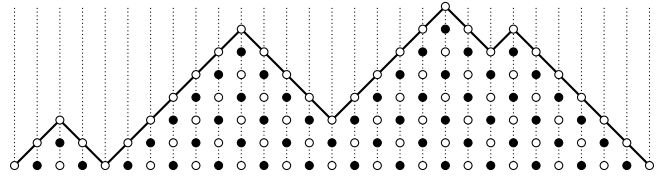
\begin{figure}[h]
	\begin{tikzpicture}[scale=0.3,rotate=0]
		\foreach\x in {0,...,28}
		{\draw[fill=white][ thin, densely dotted](\x,0)--(\x,7);}

		\draw[fill=white] [thick] (0,0) -- (1,1)--(2,2);
		
		\draw[fill=white] [thick] (2,2)--(3,1);
		\draw[fill=white] [thick] (3,1)--(4,0);
		\draw[fill=white] [thick] (4,0)--(10,6);
		\draw[fill=white] [thick] (10,6)--(11,5);
		\draw[fill=white] [thick] (11,5)--(13,3);
		\draw[fill=white] [thick ] (13,3)--(14,2);
		\draw[fill=white] [thick] (14,2)--(19,7);
		\draw[fill=white] [thick ] (19,7)--(21,5);
		\draw[fill=white] [thick] (21,5)--(22,6);
		\draw[fill=white] [thick ] (22,6)--(24,4);
		\draw[fill=white] [thick ] (24,4)--(27,1);
		\draw[fill=white] [thick ] (27,1)--(28,0);

		\foreach \x in {0}
		\foreach \y in {0}
		{
			\draw[fill=white] (\x,\y) circle (5pt);
		}
		
		\foreach \x in {1}   \foreach \y in {0}{    \filldraw (\x,\y) circle (5pt);   }
		\foreach \x in {1}   \foreach \y in {1}{ \draw[fill=white]    (\x,\y) circle (5pt);   }
		
		\foreach \x in {2}   \foreach \y in {0,2}{    \draw[fill=white] (\x,\y) circle (5pt);   }
		\foreach \x in {2}   \foreach \y in {1}{   \filldraw  (\x,\y) circle (5pt);   }
		
		\foreach \x in {3}   \foreach \y in {0}{    \filldraw (\x,\y) circle (5pt);   }
		\foreach \x in {3}   \foreach \y in {1}{   \draw[fill=white]  (\x,\y) circle (5pt);   }

		\foreach \x in {4}   \foreach \y in {0}{    \draw[fill=white] (\x,\y) circle (5pt);   }    
		\foreach \x in {5}   \foreach \y in {0}{    \filldraw (\x,\y) circle (5pt);   } 
		\foreach \x in {5}   \foreach \y in {1}{   \draw[fill=white]  (\x,\y) circle (5pt);   }    
		\foreach \x in {6}   \foreach \y in {0,2}{    \draw[fill=white] (\x,\y) circle (5pt);   } 
		\foreach \x in {6}   \foreach \y in {1}{   \filldraw  (\x,\y) circle (5pt);   }    
		\foreach \x in {7}   \foreach \y in {0,2}{    \filldraw (\x,\y) circle (5pt);   }   
		\foreach \x in {7}   \foreach \y in {1,3}{    \draw[fill=white] (\x,\y) circle (5pt);   }   
		\foreach \x in {8}   \foreach \y in {0,2,4}{    \draw[fill=white] (\x,\y) circle (5pt);   }     
		\foreach \x in {8}   \foreach \y in {1,3}{    \filldraw (\x,\y) circle (5pt);   }     
		\foreach \x in {9}   \foreach \y in {0,2,4}{    \filldraw (\x,\y) circle (5pt);   }     
		\foreach \x in {9}   \foreach \y in {1,3,5}{    \draw[fill=white] (\x,\y) circle (5pt);   }     
		\foreach \x in {10}   \foreach \y in {1,3,5}{    \filldraw (\x,\y) circle (5pt);   }     
		\foreach \x in {10}   \foreach \y in {0,2,4,6}{    \draw[fill=white] (\x,\y) circle (5pt);   }     
		\foreach \x in {11}   \foreach \y in {0,2,4}{    \filldraw (\x,\y) circle (5pt);   }     
		\foreach \x in {11}   \foreach \y in {1,3,5}{    \draw[fill=white] (\x,\y) circle (5pt);   } 
		\foreach \x in {12}   \foreach \y in {1,3}{    \filldraw (\x,\y) circle (5pt);   }     
		\foreach \x in {12}   \foreach \y in {0,2,4}{    \draw[fill=white] (\x,\y) circle (5pt);   }
		\foreach \x in {13}   \foreach \y in {0,2}{    \filldraw (\x,\y) circle (5pt);   }     
		\foreach \x in {13}   \foreach \y in {1,3}{    \draw[fill=white] (\x,\y) circle (5pt);   } 
		\foreach \x in {14}   \foreach \y in {1}{    \filldraw (\x,\y) circle (5pt);   }     
		\foreach \x in {14}   \foreach \y in {0,2}{    \draw[fill=white] (\x,\y) circle (5pt);   }     
		\foreach \x in {15}   \foreach \y in {0,2}{    \filldraw (\x,\y) circle (5pt);   }     
		\foreach \x in {15}   \foreach \y in {1,3}{    \draw[fill=white] (\x,\y) circle (5pt);   }
		\foreach \x in {16}   \foreach \y in {1,3}{    \filldraw (\x,\y) circle (5pt);   }     
		\foreach \x in {16}   \foreach \y in {0,2,4}{    \draw[fill=white] (\x,\y) circle (5pt);   } 
		\foreach \x in {17}   \foreach \y in {0,2,4}{    \filldraw (\x,\y) circle (5pt);   }     
		\foreach \x in {17}   \foreach \y in {1,3,5}{    \draw[fill=white] (\x,\y) circle (5pt);   } 
		\foreach \x in {18}   \foreach \y in {1,3,5}{    \filldraw (\x,\y) circle (5pt);   }     
		\foreach \x in {18}   \foreach \y in {0,2,4,6}{    \draw[fill=white] (\x,\y) circle (5pt);   } 
		\foreach \x in {19}   \foreach \y in {0,2,4,6}{    \filldraw (\x,\y) circle (5pt);   }     
		\foreach \x in {19}   \foreach \y in {1,3,5,7}{    \draw[fill=white] (\x,\y) circle (5pt);   } 
		\foreach \x in {20}   \foreach \y in {1,3,5}{    \filldraw (\x,\y) circle (5pt);   }     
		\foreach \x in {20}   \foreach \y in {0,2,4,6}{    \draw[fill=white] (\x,\y) circle (5pt);   }   
		\foreach \x in {21}   \foreach \y in {0,2,4}{    \filldraw (\x,\y) circle (5pt);   }     
		\foreach \x in {21}   \foreach \y in {1,3,5}{    \draw[fill=white] (\x,\y) circle (5pt);   } 
		\foreach \x in {22}   \foreach \y in {1,3,5}{    \filldraw (\x,\y) circle (5pt);   }     
		\foreach \x in {22}   \foreach \y in {0,2,4,6}{    \draw[fill=white] (\x,\y) circle (5pt);   } 
		\foreach \x in {23}   \foreach \y in {0,2,4}{    \filldraw (\x,\y) circle (5pt);   }     
		\foreach \x in {23}   \foreach \y in {1,3,5}{    \draw[fill=white] (\x,\y) circle (5pt);   } 
		
		\foreach \x in {24}   \foreach \y in {1,3}{    \filldraw (\x,\y) circle (5pt);   }     
		\foreach \x in {24}   \foreach \y in {0,2,4}{    \draw[fill=white] (\x,\y) circle (5pt);   } 
		\foreach \x in {25}   \foreach \y in {0,2}{    \filldraw (\x,\y) circle (5pt);   }     
		\foreach \x in {25}   \foreach \y in {1,3}{    \draw[fill=white] (\x,\y) circle (5pt);   } 
		
		\foreach \x in {26}   \foreach \y in {1}{    \filldraw (\x,\y) circle (5pt);   }     
		\foreach \x in {26}   \foreach \y in {0,2}{    \draw[fill=white] (\x,\y) circle (5pt);   } 
		\foreach \x in {27}   \foreach \y in {0}{    \filldraw (\x,\y) circle (5pt);   }     
		\foreach \x in {27}   \foreach \y in {1}{    \draw[fill=white] (\x,\y) circle (5pt);   } 
		
		\foreach \x in {28}   \foreach \y in {}{    \filldraw (\x,\y) circle (5pt);   }     
		\foreach \x in {28}   \foreach \y in {0}{    \draw[fill=white] (\x,\y) circle (5pt);   } 
	\end{tikzpicture}
	
	\caption{A Dyck path with black and white nodes in `checkerboard' mode.}\label{susi2}
\end{figure}

The generating functions are a bit more difficult than in the previous instance. We do not need an auxiliary sequence, and only consider $A_\ell(z)$ related to path with 
white-height $\le \ell$. It is not obvious to link it to $A_{\ell-1}(z)$. We must also consider the ordinary height (at an $x$-coordinate (=abscissa)), which is just the corresponding
$y$-value (=ordinate) of the path. A path counted by $A_\ell(z)$ will be decomposed into sojourns, namely the pieces between returns to the $x$-axis. The following drawings might clarify things a bit:

\begin{figure}[h]
 \begin{tikzpicture}[scale=0.3]
  \foreach\x in {0,...,20}
  {\draw[fill=white][ thin, densely dotted](\x,0)--(\x,4);}

  \draw[red] [thick] (0,0) -- (1,1);
  \draw[fill=white] [thick]   (1,1)--(2,2);
  
  \draw[fill=white] [thick] (2,2)--(3,1);
  \draw[fill=white] [thick] (3,1)--(6,4);
  \draw[fill=white] [thick] (6,4)--(8,2);
    \draw[fill=white] [thick] (8,2)--(10,4);
    \draw[fill=white] [thick] (10,4)--(11,3);
    \draw[fill=white] [thick] (11,3)--(12,4);
    \draw[fill=white] [thick] (12,4)--(15,1);
   \draw[fill=white] [thick] (15,1)--(17,3);
       \draw[fill=white] [thick] (17,3)--(19,1);
              \draw[red] [thick] (19,1)--(20,0);

   \draw[fill=white] (0,0) circle (5pt);
  \draw[fill=white] (1,1) circle (5pt);
    \draw[fill=white] (2,2) circle (5pt);
\draw[fill=white] (3,1) circle (5pt);
  \draw[fill=white] (4,2) circle (5pt);
    \draw[fill=white] (5,3) circle (5pt);  
    \draw[fill=white] (6,4) circle (5pt);    \draw[fill=white] (6,0) circle (5pt);
      \draw[fill=white] (7,3) circle (5pt);
        \draw[fill=white] (8,2) circle (5pt);
      \draw[fill=white] (9,3) circle (5pt);
        \draw[fill=white] (10,4) circle (5pt);    \draw[fill=white] (10,0) circle (5pt);
          \draw[fill=white] (11,3) circle (5pt);
            \draw[fill=white] (12,4) circle (5pt);    \draw[fill=white] (12,0) circle (5pt);
              \draw[fill=white] (13,3) circle (5pt);
            \draw[fill=white] (14,2) circle (5pt);
              \draw[fill=white] (15,1) circle (5pt);
                \draw[fill=white] (16,2) circle (5pt);
                  \draw[fill=white] (17,3) circle (5pt);
                    \draw[fill=white] (18,2) circle (5pt);
                  \draw[fill=white] (19,1) circle (5pt);
                    \draw[fill=white] (20,0) circle (5pt);
  
 \end{tikzpicture}\hskip 1cm%
	\begin{tikzpicture}[scale=0.3]
		\foreach\x in {1,...,19}
		{\draw[fill=white][ thin, densely dotted](\x,1)--(\x,4);}

		\draw[fill=white] [thick]  (1,1)--(2,2);
		
		\draw[fill=white] [thick] (2,2)--(3,1);
		\draw[fill=white] [thick] (3,1)--(6,4);
		\draw[fill=white] [thick] (6,4)--(8,2);
		\draw[fill=white] [thick] (8,2)--(10,4);
		\draw[fill=white] [thick] (10,4)--(11,3);
		\draw[fill=white] [thick] (11,3)--(12,4);
		\draw[fill=white] [thick] (12,4)--(15,1);
		\draw[fill=white] [thick] (15,1)--(17,3);
		\draw[fill=white] [thick] (17,3)--(19,1);

		\draw[fill=white] (1,1) circle (5pt);
		\draw[fill=white] (2,2) circle (5pt);
		\draw[fill=white] (3,1) circle (5pt);
		\draw[fill=white] (4,2) circle (5pt);
		\draw[fill=white] (5,3) circle (5pt);  
		\draw[fill=white] (6,4) circle (5pt);  		\draw[fill=black] (6,1) circle (5pt);  
		\draw[fill=white] (7,3) circle (5pt);
		\draw[fill=white] (8,2) circle (5pt);
		\draw[fill=white] (9,3) circle (5pt);
		\draw[fill=white] (10,4) circle (5pt); \draw[fill=black] (10,1) circle (5pt);  
		\draw[fill=white] (11,3) circle (5pt);
		\draw[fill=white] (12,4) circle (5pt); \draw[fill=black] (12,1) circle (5pt);  
		\draw[fill=white] (13,3) circle (5pt);
		\draw[fill=white] (14,2) circle (5pt);
		\draw[fill=white] (15,1) circle (5pt);
		\draw[fill=white] (16,2) circle (5pt);
		\draw[fill=white] (17,3) circle (5pt);
		\draw[fill=white] (18,2) circle (5pt);
		\draw[fill=white] (19,1) circle (5pt);

	\end{tikzpicture}
	
	\caption{A sojourn of height 4 and white-height 3. After chopping off first and last step, the height is 3 and  the white-height is 2}
\end{figure}
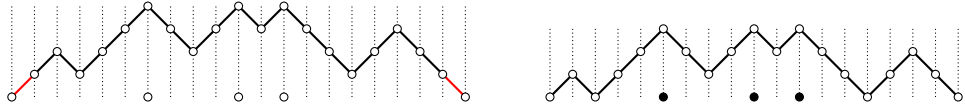


\begin{figure}[h]
	\begin{tikzpicture}[scale=0.3]
		\foreach\x in {3,...,15}
		{\draw[fill=white][ thin, densely dotted](\x,1)--(\x,4);}

				\draw[red] [thick] (3,1)--(4,2);
		\draw[fill=white] [thick] (4,2)--(6,4);
		\draw[fill=white] [thick] (6,4)--(8,2);
		\draw[fill=white] [thick] (8,2)--(10,4);
		\draw[fill=white] [thick] (10,4)--(11,3);
		\draw[fill=white] [thick] (11,3)--(12,4);
		\draw[fill=white] [thick] (12,4)--(14,2);
				\draw[fill=white,red] [thick] (14,2)--(15,1);

		\draw[fill=white] (3,1) circle (5pt);
		\draw[fill=white] (4,2) circle (5pt);
		\draw[fill=white] (5,3) circle (5pt);  
		\draw[fill=white] (6,4) circle (5pt);  		\draw[fill=black] (6,1) circle (5pt);  
		\draw[fill=white] (7,3) circle (5pt);
		\draw[fill=white] (8,2) circle (5pt);
		\draw[fill=white] (9,3) circle (5pt);
		\draw[fill=white] (10,4) circle (5pt); \draw[fill=black] (10,1) circle (5pt);  
		\draw[fill=white] (11,3) circle (5pt);
		\draw[fill=white] (12,4) circle (5pt); \draw[fill=black] (12,1) circle (5pt);  
		\draw[fill=white] (13,3) circle (5pt);
		\draw[fill=white] (14,2) circle (5pt);
		\draw[fill=white] (15,1) circle (5pt);

	\end{tikzpicture}
\hskip 1cm%
	\begin{tikzpicture}[scale=0.3]
	\foreach\x in {4,...,14}
	{\draw[fill=white][ thin, densely dotted](\x,2)--(\x,4);}

	\draw[fill=white] [thick] (4,2)--(6,4);
	\draw[fill=white] [thick] (6,4)--(8,2);
	\draw[fill=white] [thick] (8,2)--(10,4);
	\draw[fill=white] [thick] (10,4)--(11,3);
	\draw[fill=white] [thick] (11,3)--(12,4);
	\draw[fill=white] [thick] (12,4)--(14,2);

	\draw[fill=white] (4,2) circle (5pt);
	\draw[fill=white] (5,3) circle (5pt);  
	\draw[fill=white] (6,4) circle (5pt);  
	\draw[fill=white] (7,3) circle (5pt);
	\draw[fill=white] (8,2) circle (5pt);
	\draw[fill=white] (9,3) circle (5pt);
	\draw[fill=white] (10,4) circle (5pt); 
	\draw[fill=white] (11,3) circle (5pt);
	\draw[fill=white] (12,4) circle (5pt); 
	\draw[fill=white] (13,3) circle (5pt);
	\draw[fill=white] (14,2) circle (5pt);

\end{tikzpicture}
	
	\caption{A sojourn of height 3 and white-height 3. After chopping off first and last step, the height is 2 and  the white-height is still 2.}
\end{figure}
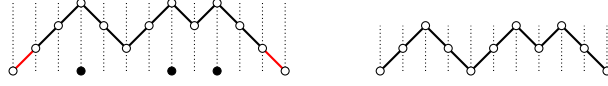

In general if the white-height of a sojourn is $=\ell$, the height of this sojourn is $=2\ell-1$ or $=2\ell-2$. By summing, if the 
white-height of a sojourn is $\le\ell$, the height of this sojourn is $\le 2\ell-1$. The quantity $H_i$ corresponding to paths of height $\le i$ is well-known.
It follows from the recursion
\begin{equation*}
	H_i=\frac{1}{1-zH_{i-1}},\ H_0=1 \quad\Longrightarrow\quad H_i= (1+u)\frac{1-u^{i+1}}{1-u^{i+2}}.
\end{equation*}
And finally
\begin{align*}
A_\ell&=\frac{1}{1-zH_{2\ell-1}}=\frac{1}{1-\dfrac{u}{(1+u)^2} (1+u)\dfrac{1-u^{2\ell}}{1-u^{2\ell+1}}}\\
&=\frac{(1+u)(1-u^{2\ell+1})}{(1+u)(1-u^{2\ell+1})-u(1-u^{2\ell})}
=(1+u)\frac{1-u^{2\ell+1}}{1-u^{2\ell+2}}.
\end{align*}
Further,
\begin{equation*}
A_\infty-A_\ell=(1+u)-(1+u)\frac{1-u^{2\ell+1}}{1-u^{2\ell+2}}=\frac{1-u^2}{u}\frac{u^{2\ell+2}}{1-u^{2\ell+2}}.
\end{equation*}
The generating function
\begin{equation*}
\frac{1-u^2}{u}\frac{u^{2\ell}}{1-u^{2\ell}}
\end{equation*}
refers to paths with white-height $\ge \ell$ in the checkerboard model. 

Pulling out coefficients can be done as before;
\begin{align*}
	[z^n]\frac{1-u^2}{u}\frac{u^{2\ell}}{1-u^{2\ell}}&=\frac1{2\pi i}\oint\frac{dz}{z^{n+1}}\frac{1-u^2}{u}\frac{u^{2\ell}}{1-u^{2\ell}}\\
	&=\frac1{2\pi i}\oint\frac{du(1-u)(1+u)^{2n+2}}{(1+u)^3u^{n+1}}\frac{1-u^2}{u}\sum_{k\ge1}u^{2k\ell}\\
	&=\sum_{k\ge1}\frac1{2\pi i}\oint\frac{du(1-u)^2(1+u)^{2n}}{u^{n+2-2k\ell}}\\
	&=\sum_{k\ge1}[u^{n+1-2k\ell}](1-2u+u^2)(1+u)^{2n}\\
	&=\sum_{k\ge1}\biggl[\binom{2n}{n+1-2k\ell}-2\binom{2n}{n-2k\ell}+\binom{2n}{n-1-2k\ell}\biggr].
\end{align*}
And
\begin{equation*}
	\sum_{\ell\ge1}\frac{e^{-2\ell t}}{1-e^{-2\ell t}}=\sum_{\ell\ge1}\sum_{k\ge1}e^{-2k\ell t};
\end{equation*}
the Mellin transform is
\begin{equation*}
\mathcal{M}\sum_{\ell\ge1}\sum_{k\ge1}e^{-2k\ell t}=\Gamma(s)2^{-s}\zeta^2(s).
\end{equation*}
Once again we must collect residues of
\begin{equation*}
	\Gamma(s)2^{-s}\zeta^2(s)t^{-s}.
\end{equation*}
Applying the inverse Mellin transform, multiplying with the factor $\frac{1-u^2}{u}$ and rewrite the expansion in the $z$-world leads to 
\begin{equation*}
-\frac12\log(1-4z)+\gamma-2\log(2)+\sqrt{1-4z}+\textsf{further terms}.
\end{equation*}
And after translating the expansion around $z\sim\frac14$ to an expansion for $n\to\infty$ (singularity analysis, transfer theorems) and finally dividing by the Catalan numbers, we get
\begin{equation*}
\textsf{Average white-height}=	\frac12\sqrt{\pi n}-\frac12+\textsf{further terms},
\end{equation*}
assuming that all Dyck paths of semi-length $n$ are equally likely. 

\section{White-height and Motzkin paths with two types of level steps}

We start by redrawing the essence of the white-height in the checkerboard sense.

\begin{figure}[h]
	\begin{tikzpicture}[scale=0.7]

		\foreach\x in {1,3,5,7,9,11,13}
		{\draw[red] [thick,-latex] (\x,1)--(\x+1,1);
			\draw[red] [thick,-latex] (\x,2)--(\x+1,2);
			\draw[red] [thick,-latex] (\x,3)--(\x+1,3);
			\draw[red] [thick,-latex] (\x,4)--(\x+1,4);
			\draw[red] [thick,-latex] (\x+1,1)--(\x+2,2);
						\draw[red] [thick,-latex] (\x+1,2)--(\x+2,3);
									\draw[red] [thick,-latex] (\x+1,3)--(\x+2,4);
			}
		
		\foreach\x in {1,3,5,7,9,11,13}
		{\draw[teal] [thick,-latex] (\x+1,2)--(\x+2,2);
			\draw[teal] [thick,-latex] (\x+1,1)--(\x+2,1);
			\draw[teal] [thick,-latex] (\x+1,3)--(\x+2,3);
			\draw[teal] [thick,-latex] (\x+1,4)--(\x+2,4);
			\draw[teal] [thick,-latex] (\x,4)--(\x+1,3);
			\draw[teal] [thick,-latex] (\x,3)--(\x+1,2);
			\draw[teal] [thick,-latex] (\x,2)--(\x+1,1);
		
		}
		
		\foreach\x in {1,...,15}
		{\draw[fill=white] (\x,1) circle (5pt);
			\draw[fill=white] (\x,2) circle (5pt);
			\draw[fill=white] (\x,3) circle (5pt);
			\draw[fill=white] (\x,4) circle (5pt);}
		
		 \draw (18,4) node[thick] {white-height \bf 4}		;
		 		 \draw (18,3) node[thick] {white-height \bf 3}		;
		 		 		 \draw (18,2) node[thick] {white-height \bf 2}		;
		 		 		 		 \draw (18,1) node[thick] {white-height \bf 1}		;
		\foreach\x in {0,...,14}
		 		 		 		 \draw (\x+1,0.2) node[] {\scriptsize{$\x$}}		;

	\end{tikzpicture}

	\caption{How red up-steps and green down-steps work.}
\end{figure}
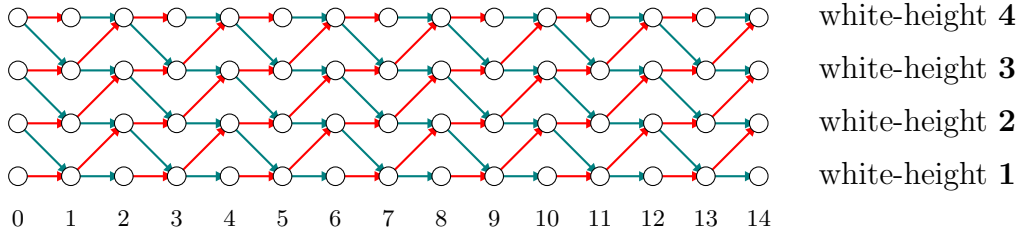
\begin{figure}[h]
	\begin{tikzpicture}[scale=0.7]

		\foreach\x in {1,3,5,7,9,11,13}
		{\draw[thick] [-latex] (\x,1)--(\x+2,1);
					\draw[red, ultra thick] [-latex] (\x,2)--(\x+2,2);
								\draw[red, ultra thick] [-latex] (\x,3)--(\x+2,3);
											\draw[red, ultra thick] [-latex] (\x,4)--(\x+2,4);
			\draw[] [thick,-latex] (\x,1)--(\x+2,2);
			\draw[] [thick,-latex] (\x,2)--(\x+2,3);
			\draw[] [thick,-latex] (\x,3)--(\x+2,4);
				\draw[] [thick,latex-] (\x+2,1)--(\x+0,2);
			\draw[] [thick,latex-] (\x+2,2)--(\x,3);
			\draw[] [thick,latex-] (\x+2,3)--(\x,4);
		}

		\foreach\x in {1,3,...,15}
		{\draw[fill=white] (\x,1) circle (5pt);
			\draw[fill=white] (\x,2) circle (5pt);
			\draw[fill=white] (\x,3) circle (5pt);
			\draw[fill=white] (\x,4) circle (5pt);}
		
		\draw (18,4) node[thick] {white-height \bf 4}		;
		\draw (18,3) node[thick] {white-height \bf 3}		;
		\draw (18,2) node[thick] {white-height \bf 2}		;
		\draw (18,1) node[thick] {white-height \bf 1}		;
		
		\foreach\x in {0,2,...,14}
		\draw (\x+1,0.2) node[] {\scriptsize{$\x$}}		;
			 
		\draw [green,thick] (0.5,1.66)--(15.5,1.66);
		
	\end{tikzpicture}

	\caption{Red stands for 2 versions of horizontal double steps, single steps are not shown.}
	\label{anika}
\end{figure}
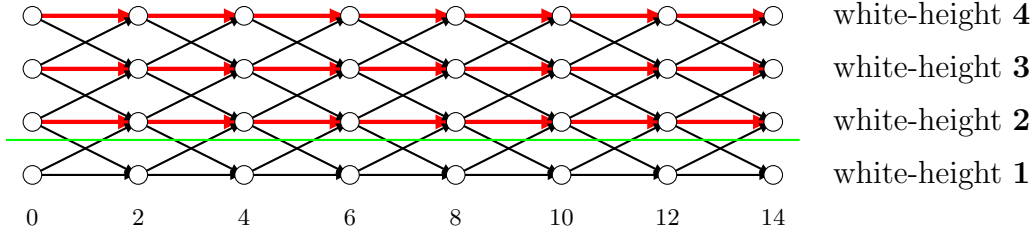

We use the variable $Z$ for  a step in the lattice \ref{anika}, which are double steps in the original lattice. First, we ignore the layer below the green line. The situation
has a Motzkin path flavour. The equation of interest is
\begin{equation*}
F=1+2ZF+Z^2F^2=(1+ZF)^2=\frac{1-\sqrt{1-4Z}}{2Z^2}-\frac1Z
\end{equation*}
Now, including the layer below the green line
\begin{equation*}
G=1+ZG+Z^2FG=\frac{1-\sqrt{1-4z}}{2z},
\end{equation*}
which is the generating function of the Dyck numbers.

This scheme can be seen in the context of the height. 
\begin{equation*}
	c_h=1+2Zc_h+Z^2c_{h-1}c_h=\frac{1}{1-2Z-Z^2c_{h-1}},\quad c_1=\frac{1}{1-2Z-Z^2c_{0}},\ c_0=0.
\end{equation*}
Then we need
\begin{equation*}
F_h=\frac{c_h}{Z^2}=(1+u)^2\frac{1-u^{2h}}{1-u^{2h+2}},
\end{equation*}
as can be seen be induction. Further
\begin{equation*}
G_{h+1}=1+ZG_{h+1}+Z^2F_{h}G_{h+1}=\frac{1}{1-Z-Z^2F_{h}}=(1+u)\frac{1-u^{2h+2}}{1-u^{2h+3}},
\end{equation*}
again by induction. To say it again:
\begin{equation*}
G_{h}=(1+u)\frac{1-u^{2h}}{1-u^{2h+1}},
\end{equation*}
which is generating function of paths with white-height $\le h$.
Here is a  special case:
\begin{equation*}
G_{1}=1+ZG_{h+1}+Z^2F_{h}G_{h+1}=\frac{1}{1-Z}=(1+u)\frac{1-u^{2}}{1-u^{3}}
\end{equation*}
since only flat steps are allowed. Further, we see that 
\begin{equation*}
G_\infty=\lim_{h\to\infty}(1+u)\frac{1-u^{2h}}{1-u^{2h+1}}=1+u=\frac{1-\sqrt{1-4z}}{2z}.
\end{equation*}

We obtained a continued fraction
\begin{equation*}
\cfrac{1}{1-Z-\cfrac{Z^2}{1-2Z-\cfrac{Z^2}{1-2Z-\cfrac{Z^2}{\dots}}}}=\frac{1-\sqrt{1-4Z}}{2Z},
\end{equation*}
which can be compared with
\begin{equation*}
	\cfrac{1}{1-\cfrac{z}{1-\cfrac{z}{1-\cfrac{z}{\dots}}}}=\frac{1-\sqrt{1-4z}}{2z}.
\end{equation*}

The combinatorial meaning of $z$ (previous section) and $Z$ (this section) is slightly different; $z$ refers to the half length, and $Z$ to the number of double-steps.
That explains the slight diffence of
\begin{equation*}
	(1+u)\frac{1-u^{2h+1}}{1-u^{2h+2}}\ \text{ versus }\  (1+u)\frac{1-u^{2h}}{1-u^{2h+1}}.
\end{equation*}
Thus, asymptotics are  already included in the first section.

We have seen in this section how 2-Motzkin paths and Dyck paths can be  be linked. Such a relation is surely known, I can think about the reference
\cite{DeutschShapiro}.

\end{document}